\documentclass[12pt]{article}

\usepackage{latexsym,palatino,amsthm,amssymb,amsmath,amsfonts,graphicx,color,tikz}
\usepackage[colorlinks=true,citecolor=black,linkcolor=black,urlcolor=gray]{hyperref}

\allowdisplaybreaks

\newcommand{\beq}{\begin{equation}}
	\newcommand{\eeq}{\end{equation}}
\newcommand{\bqa}{\begin{eqnarray}}
	\newcommand{\eqa}{\end{eqnarray}}
\newcommand{\nn}{\nonumber}

\begin{document}
	
	\section*{Matrices with integer eigenvalues for all permutations of coefficients (thanks to Pythagoras!)}
	\vspace{-.30cm}
	
	\begin{center}
		\textbf{Michael J.W. Hall}\footnote{%
		Michael Hall is an honorary senior lecturer and keen puzzlist in the Department of Theoretical Physics, in the Research School of Physics at the Australian National University --- see also \url{https://www.researchgate.net/profile/Michael-Hall-79}.}
	\end{center}
	
	\medskip

	\section*{Introduction}
	
	In Year 10 of high school my maths teacher, Miss Thornley,  began each class by reeling off ten problems made up on the spot for us to quickly answer --- such as factoring $x^2-5x+6$, or solving $x-3y=7$ and $2x+y=0$ for $x$ and $y$. She would then similarly reel off the solutions --- which, to make things easier for us, were always integers --- while we scored our neighbours' responses. This exercise not only kicked our brains into gear, but was impressive in the way she generated and remembered the questions.
	
	Miss Thornley didn't teach us about matrices --- arrays of numbers such as 
	$\begin{pmatrix}
		a & b \\ c & d
	\end{pmatrix}$, 
	which represent transformations that send straight lines to straight lines,
	including rotations, reflections and stretchings --- but if she had, I expect she would have asked the class to find the scaling factors by which a given matrix  can transform straight lines. These factors are determined by calculating the ``eigenvalues'' of the matrix, and hence a method for  quickly generating matrices with both integer coefficients and integer eigenvalues would be useful in this and in other teaching contexts.
	
	Methods of this sort are in fact known, which generate suitable matrices in terms of some set of parameters~\cite{cronin,towse,phan}. An alternative method is given here, however, which is more efficient in the sense that each choice of parameter values generates more than one suitable matrix. Indeed, {\it every} permutation of the numbers in a matrix generated by this method yields a new matrix that also has integer eigenvalues. In this way, for example, a maths teacher can give 24 different $2\times2$ matrices to 24 respective students, all having the same 4 coefficients $a,b,c,d$ in some order, and all having integer eigenvalues. 
	
	Further, this alternative method generates such matrices simply and directly from Pythagorean triples (or triads), i.e., from the sides of right--angled triangles having integer lengths. For example, each and every permutation of the numbers appearing in the matrix $\begin{pmatrix} 12 & 6 \\ 7 & 1 \end{pmatrix}$, such as $\begin{pmatrix} 6 &1 \\ 7 & 12 \end{pmatrix}$, has integer eigenvalues generated by the Pythagorean triple $(5,12,13)$. Note the construction in~\cite{phan} also uses Pythagorean triples, but does not have the permutation property (as demonstrated, e.g., by the simple choice $a=0, n=1$ in Theorem~1 of~\cite{phan}). The method is derived below, followed by some explorations of the results. No prior knowledge of matrix theory is required.

	\section*{The basic equations}
	
	The only result needed from matrix theory is that the eigenvalue equation for a $2\times 2$ matrix 
	$\begin{pmatrix}
		a & b \\ c & d
	\end{pmatrix}$
	is the quadratic equation $\lambda^2-(a+d)\lambda+(ad-bc)=0$, which has solutions
	\beq \nn
	\lambda_\pm=\frac{a+d \pm \sqrt{(a+d)^2-4(ad-bc)}}{2}
	\eeq
	for the corresponding eigenvalues $\lambda_+, \lambda_-$. Note for integer coefficients $a,b,c,d$ that the square--root term is rational if and only if it is a whole number $u$~\cite{inf}, and that in this case $u$ has the same parity (odd or even) as  $a+d$, so that $a+d\pm u$ is even. It immediately follows that the eigenvalues are also integers if and only if the condition
	\beq \label{c1}
	u^2 = (a+d)^2-4(ad-bc) = (a-d)^2+ 4bc
	\eeq
	is satisfied for some integer $u$. 
	
	The above condition is invariant under the interchanges $a\leftrightarrow d$ and/or $b\leftrightarrow c$, i.e., under swapping of the diagonal and/or off--diagonal coefficients of the matrix. Hence, Eq.~(\ref{c1}) applies to 4 possible permutations of the coefficients. More generally, {\it all} 24 matrices that can be obtained by permutations of $a,b,c$ and $d$ will have integer eigenvalues if and only if Eq.~(\ref{c1}) and its 5 permutations 
	\begin{align}
		v^2 &= (b-c)^2 + 4ad      \label{c2} \\
		w^2 &=  (a-b)^2 + 4cd     \label{c3} \\
		x^2 &=   (c-d)^2 + 4ab  \label{c4}\\
		y^2 &= (a-c)^2 + 4bd \label{c5}\\
		z^2 &= (b-d)^2 + 4ac \label{c6}
	\end{align}
	are satisfied for some choice of integers $u,v,w,x,y,z$.
	
	The aim is therefore to find sets of integers $\{a,b,c,d\}$ for which conditions~(\ref{c1})--(\ref{c6}) hold. It will be seen that a large class of solutions are generated by Pythagorean triples. A second class of solutions is noted at the end of this article.
	
	\section*{Solving the equations via Pythagorean triples}
	
	The promised ``Pythagorean'' solutions are obtained by making a simplifying assumption or ansatz:
	\beq \label{ansatz}
	a+d=b+c =t
	\eeq
	for some integer $t$. That is, the sum of two of the coefficients is equal to the sum of the other two coefficients. 
	
	Under this ansatz we have $a-b=c-d$ and $a-c=b-d$, so that conditions~(\ref{c3})--(\ref{c6}) are automatically satisfied, with
	\beq \label{eig1}
	w=c+d,~~ x=a+b,~~ y=b+ d,~~ z= a+c .
	\eeq
	Further, condition~(\ref{c1}) simplifies to
	\begin{align}
		u^2 &= (a+d)^2-4(ad-bc) \nn\\
		&= (a+d)^2-4[ad-b(a-b+d)] \nn\\
		&= t^2 + 4(a-b)(b-d) ,
		\label{y}
	\end{align}
	while condition~(\ref{c2}) similarly simplifies to
	\beq \label{z}
	v^2 = t^2 - 4(a-b)(b-d).
	\eeq
	Hence, the problem is reduced to finding suitable solutions of Eqs.~(\ref{y}) and~(\ref{z}).
	
	First, subtracting Eq.~(\ref{z}) from Eq.~(\ref{y}) gives
	\beq \nn
	u^2-v^2 = (u+v)(u-v) = 8(a-b)(b-d) .
	\eeq
	Noting  $u+v$ and $u-v$ have the same parity (odd or even) and an even product, they must both be even. Hence $u+v=2r$ and $u-v=2s$ for two integers $r$ and $s$. It follows that
	\beq \label{yz} 
	u=r+s,~~ v=r-s, ~~ rs = 2(a-b)(b-d).
	\eeq
	Further, substitution into the sum of Eqs.~(\ref{y}) and~(\ref{z})  yields $u^2+v^2=2t^2=(r+s)^2 + (r-s)^2=2r^2+2s^2$, i.e.,
	\beq \label{py}
	r^2+s^2 = t^2 .
	\eeq
	Thus $(r,s,t)$ forms a Pythagorean triple. 
	
	Defining $k=b-d$ and $l=a-b$, then $k+l=a-d=2a-t$ and $k-l=2b-a-d=2b-t$ using ansatz~(\ref{ansatz}), while $rs=2kl$ from Eq.~(\ref{yz}). The 4 matrix coefficients therefore have the form
	\begin{align} 
		\{ a,b,c,d\} &= \left\{ \frac{t+ k+l }{2}, \frac{t+ k- l}{2}, \frac{t- k+ l}{2}, \frac{t- k- l}{2} \right\} \nn\\
		\label{sol}
		&= \left\{ \frac{t\pm k\pm l}{2} \right\},~~~{\rm with}~~kl=\frac{rs}{2}.
	\end{align}
	
	It remains to be ensured, for a given Pythagorean triple $(r,s,t)$, that factors $k$ and $l$ of $\frac{rs}{2}$ exist such that $a,b,c,d$ are indeed integers. There are in fact many ways of choosing such $k$ and $l$, as considered further below. However, to give an immediate example, note that for any Pythagorean triple $(r,s,t)$ it is always the case that at least one of $r$ or $s$ has the same parity as $t$ and the other is divisible by 4.\footnote{ This follows immediately from the representation of Pythagorean triples in Eq.~(\ref{rep}) below. More directly, note that the possibility of a triple with $r$ and $s$ both odd and $t$ even  is ruled out by implying $0=r^2+s^2-t^2=4p+2$ for some integer $p$. This leaves two possibilities: $r$ and $s$ are either both even or have opposite parities. In the first case $t$ is also even, and so one can keep dividing the triple by 2 until the second possibility is realised. Finally, in the  latter case suppose that $r=2p+1$ is odd, without loss of generality, so that $t=2q+1$ is also odd. Hence $s^2=t^2-r^2=8(T_p-T_q)$, where $T_m$ denotes the triangular number $\frac12 m(m+1)$. Thus $s^2$ is divisible by 8, which implies $s$ is divisible by 4 as claimed.} 
	Hence one can always guarantee integer coefficients in Eq.~(\ref{sol}) by choosing $k=r$ to have the same parity as $t$ and $l=\frac{s}{2}$, to give the simple {\it canonical} solution 
	\beq \label{can}
	\{ a,b,c,d\} =  \left\{ \frac{t\pm r\pm \frac{s}{2}}{2} \right\}
	\eeq
	for any Pythagorean triple $(r,s,t)$.
	
	Thus, for example, the triples $(3,4,5)$ and $(5,12,13)$ generate the respective canonical solutions 
	\beq \label{ex}
	\{ a,b,c,d\} =  \left\{ 5,3,2,0 \right\},~~~\{ a,b,c,d\} =\left\{ 12,6,7,1 \right\} ,
	\eeq 
	Any $2\times2$ matrix formed by some permutation of these coefficients will have integer eigenvalues.  Other examples are discussed below.
	
	\section*{Further properties and examples}
	
	First, a little more about the canonical solution in Eq.~(\ref{can}).
	In particular, from the unique representation
	\beq \label{rep}
	r=(f^2-g^2)h, ~~s=2fgh, ~~t=(f^2+g^2)h
	\eeq
	of Pythagorean triples, where $f$ and $g$ are coprime and have opposite parities, and $h$ is an arbitrary common factor~\cite{triple}, the canonical solution reduces to  $\{ a,b,c,d\} = h\{ f^2 \pm \frac12 fg, g^2\pm \frac12 fg\}$. Hence, taking $f=2m$ without loss of generality and rewriting $g$ as $n$, the canonical solutions have the simple alternative form
	\beq \label{altcan}
	\{ a,b,c,d\} = \{ 4m^2 \pm mn, n^2\pm mn \} ,
	\eeq
	up to an arbitrary common factor, for any coprime integers  $m$ and $n$ with $n$ odd.
	
	The above alternative form is particularly easy to remember and use, and recovers the examples in Eq.~(\ref{ex}) for the respective choices $m=n=1$ and $m=1, n=3$. Note that the coefficients in the second example are strictly positive. From Eq.~(\ref{altcan}) one has more generally that
	\beq
	a,b,c,d>0 {\rm ~~ if ~and ~only ~if~~} 4|m|>|n|>|m| .
	\eeq
	
	There are also many noncanonical solutions corresponding to Eq.~(\ref{sol}). For example, whenever $t$ is odd one can choose $k=\frac{rs}{2}$ and $l=1$, to give the corresponding solution
	\beq \label{altsol}
	\{ a,b,c,d\} = \left\{ \frac{t\pm  \frac{rs}{2} \pm 1}{2} \right\}
	\eeq
	(and any multiple thereof), for any Pythagorean triple $(r,s,t)$ with $t$ odd.
	For the triples $(3,4,5)$ and $(5,12,13)$, for example, this  generates the respective solutions
	\beq \label{ex2}
	\{ a,b,c,d\} =  \left\{ 6, 5, 0, -1 \right\},~~~\{ a,b,c,d\} =\left\{22, 21, -8, -9 \right\} ,
	\eeq 
	which may be compared to the canonical examples in Eq.~(\ref{ex}). Again, any $2\times2$ matrix formed by some permutation of these coefficients will have integer eigenvalues.
	
	Finally, even for non--integer rational choices of $k$ and $l$, for which Eq.~(\ref{sol}) generates rational rather than integer values for $a,b,c,d$, one can nevertheless multiply these values by a suitable common factor to obtain integer solutions. For example, for the Pythagorean triple $(3,4,5)$, the choices $k=4, l=\frac{3}{2}$ and $k=\frac{9}{2}, l=\frac{4}{3}$ give fractional values for $\{a,b,c,d\}$ in Eq.~(\ref{sol}), which can be multiplied by 4 and 12 respectively to give the corresponding integer solutions
	\beq \label{ex3}
	\{ a,b,c,d\} =  \left\{ 21, 15, 5, -1 \right\},~~~\{ a,b,c,d\} =\left\{65, 49, 11, -5 \right\} .
	\eeq 
	More generally, for any triple $(r,s,t)$  the choice $k=\frac{p}{q}, l=\frac{qrs}{2p}$ and multiplying factor $pq$ generates the solution
	\beq \label{count}
	\{ a,b,c,d\} = \left\{ pqt\pm p^2 \pm \frac{q^2rs}{2} \right\}
	\eeq
	for any rational number $\frac{p}{q}$ (where the corresponding eigenvalues in Eq.~(\ref{dadah}) below must also be multiplied by $pq$).
	In this manner each Pythagorean triple can in fact generate a countable infinity of nontrivial solutions, although the canonical solutions in Eqs.~(\ref{can}) and~(\ref{altcan}) remain the simplest.

	\section*{And the eigenvalues are \dots}
	
	It is,  of course, of interest to explicitly determine the integer eigenvalues generated by a given Pythagorean triple $(r,s,t)$. This is straightforward even for the general solution in Eq.~(\ref{sol}).
	
	Recall first that the eigenvalues of the $2\times 2$ matrix 
	$\begin{pmatrix}
		a & b \\ c & d
	\end{pmatrix}$
are given by 
	\beq \label{eig}
	\lambda_\pm=\frac{a+d \pm \sqrt{(a+d)^2-4(ad-bc)}}{2} .
	\eeq
	This formula allows us to determine the forms of the integer eigenvalues corresponding to each of the matrices
	\beq \label{perms}
	\begin{pmatrix}
		a & b \\ c & d
	\end{pmatrix} ,~~
	\begin{pmatrix}
		b & a \\ d & c
	\end{pmatrix},~~
	\begin{pmatrix}
		a & d \\ c & b
	\end{pmatrix},~~
	\begin{pmatrix}
		d & a \\ b & c
	\end{pmatrix} ,~~
	\begin{pmatrix}
		a & b \\ d & c
	\end{pmatrix},~~
	\begin{pmatrix}
		b & a \\ c & d
	\end{pmatrix},
	\eeq
	where
	\beq
	a=\frac{t+k+l}{2}, ~~~	b=\frac{t+k-l}{2}, ~~~	c=\frac{t-k+l}{2}, ~~~	d=\frac{t-k-l}{2}
	\eeq
	as per Eq.~(\ref{sol}), with $kl=\frac{rs}{2}$. In particular, for each of the matrices in Eq.~(\ref{perms})
	the corresponding square--root term in eigenvalue equation~(\ref{eig}) is given by $u,v,w,x,y$ and $z$ in  conditions~(\ref{c1})--(\ref{c6}) respectively. The corresponding pairs of eigenvalues then follow easily via the expressions for $u,v,w,x,y,z$ in Eqs.~(\ref{eig1}) and~(\ref{yz}), as
	\beq \label{dadah} 
	\left\{\frac{t\pm(r+s)}{2}\right\},~~ \left\{\frac{t\pm(r-s)}{2}\right\},~~ \{t,k\},~~\{t,-k\},~~ \{t, l\},~~ \{t, -l\}, 
	\eeq
	respectively. For the case of the canonical solution in Eq.~(\ref{can}) one has $k=r$, $l=\frac{s}{2}$.
	
	It is interesting to note that the eigenvalues of the first two matrices are independent of the choice of $k$ and $l$. Note also for each one of the matrices in Eq.~(\ref{perms}) that the further 3 matrices generated by swapping the diagonal and/or off-diagonal coefficients have the same eigenvalues, since the form of Eq.~(\ref{eig}) is invariant under such swaps.

	\section*{Conclusion}
	
	It is seen that Pythagorean triples can be used to generate matrices that have integer eigenvalues for all  permutations of their coefficients, via simple formulas such as Eqs.~(\ref{sol}), (\ref{can}) and~(\ref{altcan}). Further, while integer multiples of such matrices trivially generate additional examples, it has been also shown that just a single Pythagorean triple can generate a countable infinity of nontrivially related matrices, as per Eq.~(\ref{count}).
	
	It is natural to ask whether the ansatz $a+d=b+c$ in Eq.~(\ref{ansatz}) is a necessary as well as a sufficient condition for solutions to the basic equations~(\ref{c1})--(\ref{c6}) to exist. The answer is provided by considering the alternative ansatz $c=d=0$ in these equations, which is left to the interested reader.   
	
	Finally, it would be of interest to either find or rule out solutions for the case of $3\times3$ matrices, in addition to the trivial solution with all 9 coefficients equal.

	\vspace{0.5cm}
	
	\flushleft{Address for correspondence:\\
		\it Michael J W Hall,\\ 
		Theoretical Physics,
		Research School of Physics,\\
		Australian National University,
		Canberra, ACT 0200, Australia}

	\end{document}